\input amstex\documentstyle{amsppt}  
\pagewidth{12.5cm}\pageheight{19cm}\magnification\magstep1
\topmatter
\title Unipotent character sheaves and strata of a reductive group
\endtitle
\author G. Lusztig\endauthor
\address{Department of Mathematics, M.I.T., Cambridge, MA 02139}\endaddress
\thanks{Supported by NSF grant DMS-2153741}\endthanks
\abstract{ Let $H$ be a connected reductive group over an algebraically closed field. We define a surjective map from the set $CS(H)$ of unipotent character sheaves on $H$ (up to isomorphism) to the set of strata of
$H$. To do this we use the generalized Springer correspondence. We also give a new parametrization of $CS(H)$ in terms of data coming from bad
characteristic.} \endabstract
\endtopmatter   
\document

\define\bx{\boxed}

\define\Irr{\text{\rm Irr}}

\define\si{\sim}

\define\sqc{\sqcup}

\define\lb{\linebreak}

\define\op{\oplus}
   
\define\part{\partial}
\define\emp{\emptyset}

\define\n{\notin}

\define\m{\mapsto}
\define\do{\dots}

\define\lra{\leftrightarrow}

\define\sub{\subset}    

\define\T{\times}
\define\ti{\tilde}
\define\nl{\newline}
\redefine\i{^{-1}}

\define\un{\underline}
\define\ov{\overline}

\define\bbq{\bar{\QQ}_l}

\define\Hom{\text{\rm Hom}}
\define\End{\text{\rm End}}

\define\ind{\text{\rm ind}}

\redefine\c{\chi}
\define\g{\gamma}
\redefine\d{\delta}
\define\e{\epsilon}

\define\ph{\phi}

\define\r{\rho}
\define\s{\sigma}
\redefine\t{\tau}
\define\th{\theta}

\redefine\G{\Gamma}
\redefine\D{\Delta}

\define\Ph{\Phi}

\define\kk{\bold k}

\define\CC{\bold C}

\define\NN{\bold N}

\define\QQ{\bold Q}

\define\SS{\bold S}

\define\ca{\Cal A}

\define\cc{\Cal C}

\define\cl{\Cal L}

\define\cu{\Cal U}

\define\cw{\Cal W}

\define\cx{\Cal X}

\define\sha{\sharp}

\head Introduction\endhead
\subhead 0.1\endsubhead
Let $H$ be a reductive connected group over $\CC$.
Let $Pr=\{2,3,5,\do\}$ be the set of prime numbers;
let $\ov{Pr}=Pr\cup\{0\}$.
For $r\in Pr$ let 
$\kk_r$ be an algebraic closure of a finite field with $r$ elements
and let $H_r$ be a reductive connected group over $\kk_r$ of the same
type as $H$ and with the
same Weyl group $W$ (with set of simple reflections $\{s_i;i\in I\}$).
We set $\kk_0=\CC,H_0=H$. Let $K_r=\CC$ (if $r=0$) and $K_r=\bbq$ 
where $l\in Pr-\{r\}$ (if $r\in Pr$).
For $r\in\ov{Pr}$ 
let $CS(H_r)$ be the (finite) set of isomorphism classes of unipotent
character sheaves on $H_r$. These are certain simple perverse
$K_r$-sheaves on $H_r$, see \cite{L85}. It is known that $CS(H_r)$
is independent of $r$ in a canonical way.
Let $Str(H_r)$ be the (finite) set of strata of $H_r$, see
\cite{L15}; these are certain subsets of $H_r$ (unions of conjugacy
classes of fixed dimension) which form a partition of $H_r$. (These
subsets are locally closed in $H_r$, see \cite{C20}.)

In this paper we shall define for any $r\in\ov{Pr}$ a surjective
map
$$\t:CS(H_r)@>>>Str(H_r).\tag a$$

In the remainder of this paper (except in 1.1 and 3.3)
we shall assume that
either $H$ is quasi-simple, that is, $H$ modulo its centre is simple,
or that $H$ is a torus.
(The general case can be reduced in an obvious way to this case.)

Our definition of the map (a) is based on the generalized Springer
correspondence of \cite{L84}, especially in bad characteristic.

In 1.11 we use the map (a) to give a new parametrization of
$CS(H_r)$ which differs from the known classification \cite{L86}
in terms of two sided cells in the Weyl group.
This involves associating to each stratum a very small
collection of finite groups which come from unipotent classes in bad
characteristic. (It would be interesting to give a definition of these
finite groups and of the resulting parametrization which is purely
in characteristic $0$.)

This can be also viewed as a parametrization of 

(b) $Un(H_r(F_q))$, the set of unipotent representations of the
group $H_r(F_q)$ of $F_q$-rational points of a split form
of $H_r$ over a finite subfield $F_q$ of $\kk_r$ (with
$r\in Pr$ and $q$ a power of $r$).

Indeed, it is known that

(c) $CS(H_r),Un(H_r(F_q))$ are in a natural $1-1$ correspondence.

We will show elsewhere that similar results hold when $H$ is
replaced by a connected component of a disconnected reductive group
with identity component $H$.

A number of results of this  paper rely on the wonderful paper
\cite{S85} of Nicolas Spaltenstein and also on \cite{LS85}.

\subhead 0.2. Notation \endsubhead
Assume that $G$ is a connected reductive group over an algebraically
closed field. We denote 
by $Z_G$ the centre of $G$ and by $Z_G^0$ its identity component.
Let $G_{ad}$ be the adjoint group of $G$.

For $g\in G$ we denote by $Z_G(g)$ the centralizer of $g$ in $G$
and by $Z^0_G(g)$ its identity component. 

If $G'$ is a subgroup of $G$, we denote by $N_G(G')$ the normalizer of
$G'$ in $G$.

If $\cw$ is a Weyl group we denote by $\Irr(\cw)$ the set of
isomorphism classes of irreducible representations of $\cw$ over $\QQ$.

\head 1. Definition of the map $\t$\endhead
\subhead 1.1\endsubhead
For $r\in\ov{Pr}$ let $CS^\emp(H_r)$ be the subset of $CS(H_r)$
consisting of unipotent cuspidal character sheaves.

Let $A\in CS^\emp(H_r)$. The support of $A$ is the closure in $H_r$ of a
single orbit of $Z^0_{H_r}\T H_r$ acting on $H_r$ by
$(z,g):g_1\m zgg_1g\i$; this orbit is denoted by $\s_A$.
Let $\d(A)$ be the dimension of the variety of Borel subgroups of $H_r$
that contain a fixed element $h\in\s_A$ (this is independent of the choice of $h$). We have
$$CS^\emp(H_r)=\sqc_{d\in\NN}CS^\emp_d(H_r)$$
where $CS^\emp_d(H_r)=\{A\in CS^\emp(H_r);\d(A)=d\}$.

\proclaim{Lemma 1.2} For any $d\in\NN$ the function
$r\m\sha(CS^\emp_d(H_r))$ from $\ov{Pr}$ to $\NN$ is constant; its
value is denoted by $N_d(H)\in\NN$.
\endproclaim  
This can be deduced from the results in \cite{L86}. (When $H_{ad}$ is of
type $E_8$ or $F_4$, the results in {\it loc.cit.} are proved only under
the assumption that $r$ is not a bad prime for $H$. But the same proof
works without this assumption, by making use of \cite{S85, p.336,337}.)

\subhead 1.3\endsubhead
For $r\in\ov{Pr}$ and $J\sub I$ we fix a Levi subgroup $L_{J,r}$ of
a parabolic subgroup of $H_r$ of type $J$. (For example, $L_{I,r}=H_r$
and $L_{\emp,r}$ is a maximal torus.) We say that $J$ is {\it cuspidal}
if for some (or equivalently any) $r\in\ov{Pr}$ we have
$CS^\emp(L_{J,r})\ne\emp$. In this case $L_{J,r}$ is quasi-simple or a
torus and $J$ is uniquely determined by the type of $(L_{J,r})_{ad}$.
(This follows from the classification of cuspidal character sheaves.)
Let $W_J$ be the Weyl group of $L_{J,r}$, viewed as a parabolic subgroup
of $W$.

Let us now fix a cuspidal $J$ and $A'\in CS^\emp(L_{J,r})$.
The induced object $\ind(A')$ is a well defined semisimple perverse
sheaf on $H_r$ (see \cite{L85, \S4}); it is in fact a direct sum of
character sheaves on $H_r$. By arguments in \cite{L84, \S3,\S4},
$\End(\ind(A'))$ has a canonical decomposition as a direct
sum of lines $\op_w\cl_w$ with $w$ running through
$N_{H_r}(L_{J,r})/L_{J,r}=N_W(W_J)/W_J$ such that $\cl_w\cl_{w'}=\cl_{ww'}$ for any
$w,w'$ in $N_W(W_J)/W_J$.
One can verify that there is a unique
$A\in CS(H_r)$ such that $A$ is a summand with multiplicity one
of $\ind(A')$ and the value of the $a$-function of $W$
on the two-sided cell of $W$ attached to $A$ is equal to the value of the
$a$-function of $W_J$ on the two-sided cell of $W_J$ attached to $A'$.
Now the summand $A$ of $\ind(A')$ is stable under each $\cl_w$ and we
can choose uniquely a nonzero vector $t_w\in\cl_w$ which acts on $A$
as identity. We have $t_wt_{w'}=t_{ww'}$ for any $w,w'$ in
$N_W(W_J)/W_J$.
We see that $\End(\ind(A'))$ is canonically the group algebra of
$N_W(W_J)/W_J$ (which is known to be a Weyl group). For any
$E'\in\Irr(N_W(W_J)/W_J)$ let $A'[E']$ be the perverse sheaf
$\Hom_{N_W(W_J)/W_J}(E',\ind(A'))$ on $H_r$. This is an object of
$CS(H_r)$.

\subhead 1.4\endsubhead
Let $CS'(H_r)$ be the set of triples $(J,E',A')$ where
$J$ is a cuspidal subset of $I$, $E'\in\Irr(N_W(W_J)/W_J)$ and
$A'\in CS^\emp(L_{J,r})$. We have a bijection
$$CS'(H_r)@>\si>>CS(H_r)\tag a$$
given by $(J,E',A')\m A'[E']$.

\subhead 1.5\endsubhead
Let $r\in\ov{Pr}$. Let $\cu(H_r)$ be the set of unipotent
classes in $H_r$; for $\g\in\cu(H_r)$ the Springer correspondence
(defined for any $r$ in \cite{L84})
associates to $\g$ and the constant local system $K_r$ on $\g$
an element $e_r(\g)\in\Irr(W)$. Thus we have a well defined (injective)
map $e_r:\cu(H_r)@>>>\Irr(W)$, whose image is denoted by $\Irr_r(W)$.

Let $CS^\emp(H_r)^{un}$ be the set of all $A\in CS^\emp(H_r)$ such that 
$\s_A=Z_{H_r}^0\g_A$ where $\g_A\in\cu(H_r)$.
Let 
$$CS'(H_r)^{un}=\{(J,E',A')\in CS'(H_r);A'\in CS^\emp(L_{J,r})^{un}\}.$$
We define a map

(a) $\ti e_r:CS'(H_r)^{un}@>>>\Irr_r(W)$
\nl
as follows.
Let $(J,E',A')\in CS'(H_r)^{un}$. Then the unipotent class $\g_{A'}$
of $L_{J,r}$ is defined; the restriction of $A'$ to $\g_{A'}$ is (up to a shift) a cuspidal
local system. Now the generalized Springer correspondence \cite{L84}
associates to this cuspidal local system and to $E'$ a unipotent class
$\g$ of $H_r$ and an irreducible local system on it. By definition,
we have $\ti e_r(J,E',A')=e_r(\g)$.

\subhead 1.6\endsubhead
Let
$$\Irr_*(W)=\cup_{r\in Pr}\Irr_r(W)=\cup_{r\in\ov{Pr}}\Irr_r(W).$$
Let $r\in\ov{Pr}$. In \cite{L15} a bijection
$$Str(H_r)@>>>\Irr_*(W)\tag a$$
is defined. Using this and 1.4(a), we see that defining $\t$ in
0.1(a) is the same as defining a map
$$\un\t_r:CS'(H_r)@>>>\Irr_*(W).$$

\proclaim{Lemma 1.7} Let $d\in\NN$
be such that $N_d(H)>0$ (see 1.2).
Let $X=\{r\in Pr;CS^\emp_d(H_r)\sub CS^\emp(H_r)^{un}\}$.
Then one of the following holds.

(i) $X$ consists of a single element $r_0$.

(ii) $X=Pr$ and $d\ge1$. (In this case $H_{ad}$ is of type
$E_8,F_4$ or $G_2$, $d$ is $16,4,1$ respectively and $N_d(H)=1$.)

(iii) $X=Pr$ and $d=0$. (In this case $H$ is a torus.)

(iv) $X=\emp$. (In this case $d=0$ and $H_{ad}$ is of type $E_8$,
$F_4$ or $G_2$.)
\endproclaim
This follows from 3.2.

\subhead 1.8\endsubhead
Let $r\in\ov{Pr}$. We will now define the map $\un\t_r:CS'(H_r)@>>>\Irr_*(W)$.
In the case where $I=\emp$, this map is the bijection between
two sets with one element. Assume now that $I\ne\emp$. Let
$(J,E',A')\in CS'(H_r)$. We want to define $\un\t_r(J,E',A')$. 

Let $X$ be as in Lemma 1.7 for $L_{J,r}$ instead of $H_r$ and for
$d\in\NN$ defined by $A'\in CS^\emp_d(L_{J,r})$.

Assume first that $J\ne I,J\ne\emp$. Then $X$ is not as in
1.7(ii),(iii),(iv), hence it is as in 1.7(i). Let $r_0\in Pr$
be such that $X=\{r_0\}$.
We set $\un\t_r(J,E',A')=\ti e_{r_0}(J,E',A')$.

Next we assume that $J=I$ and $X$ is as in 1.7(i).
Let $r_0\in Pr$ be such that $X=\{r_0\}$.
We set $\un\t_r(J,E',A')=\ti e_{r_0}(J,E',A')$.

Next we assume that $J=I$ and $d,X$ are as in 1.7(ii). We have
$E'=1$. We set $\un\t_r(I,1,A')=\ti e_{r'}(I,1,A')$ where
$r'\in\ov{Pr}$ (this is independent of $r'$ by
results in \cite {S85}). 

If $J=I$ then $d$ cannot be as in 1.7(iii) since this would imply
$I=\emp$ contrary to our assumption.

Assume now that $J=I$ and $X$ is as in 1.7(iv).
We have $E'=1$. We set $\un\t_r(I,1,A')=
\text{ unit representation}$.

Finally, assume that $J=\emp$. Then $A'$ is the constant sheaf $K_r$.
For any $r'\in Pr$ we set $\ti e_{r'}(\emp,E',K_r)=E'_{r'}$.
If $E'_{r'}$ is independent of $r'$, then
$\un\t_r(\emp,E',K_r)$ is defined to be this constant value of $E'_{r'}$.
If $E'_{r'}$ is not independent of $r'$, then there is a unique
$r_0\in Pr$ such that $E'_{r'}$ is constant for $r'\in Pr-\{r_0\}$.
(This is an issue only in exceptional types where it can be checked
from the tables in \cite{S85}.) We then set
$\un\t_r(\emp,E',K_r)=E'_{r_0}$. This completes the definition of
$\un\t_r$ hence also that of $\t$ in 0.1(a).

\subhead 1.9\endsubhead
Let $r\in\ov{Pr}$. If $H_{ad}$ is of classical type or of type $E_6,E_7$ or $F_4$,
then for any $E'\in\Irr(W)$, $E'_{r'}$ is constant for $r'\in Pr-\{2\}$.
It follows that
$$\un\t_r(\emp,E',K_r)=E'_2.$$
Hence if $E'\in\Irr_2(W)$, then $\un\t_r(\emp,E',K_r)=E'.$
If $H_{ad}$ is of type $G_2$, then for any $E'\in\Irr(W)$,
$E'_{r'}$ is constant for $r'\in Pr-\{3\}$. It follows that
$$\un\t_r(\emp,E',K_r)=E'_3.$$
Hence if $E'\in\Irr_3(W)$, then $\un\t_r(\emp,E',K_r)=E'.$
We see that if $H_{ad}$ is not of type $E_8$, then $\un\t_r(\emp,E',K_r)=E'$
for $E'\in\Irr_*(W)$.
The same holds if $H_{ad}$ is of type $E_8$ (we use the tables in
\cite{S85}).

Note that $\Irr_*(W)$ can be viewed as a subset of $CS'(H_r)$ by
$E'\m(\emp,E',K_r)$. The results above show that $\un\t_r$ can be viewed
as a retraction of $CS'(H_r)$ onto its subset $\Irr_*(W)$. In
particular, $\un\t_r$ is surjective.

\subhead 1.10\endsubhead
Let $E\in\Irr_*(W)$. Let $\ov{Pr}(E)=\{r'\in\ov{Pr};E\in\Irr_{r'}(W)\}$.
For $r'\in\ov{Pr}(E)$ we denote by $\g_E$ the unique element of
$\cu_{H_{r'}}$ such that $e_{r'}(\g_E)=E$ (see 1.5). We set
$\ca_{r',E}=Z_{(H_{r'})_{ad}}(u)/Z_{(H_{r'})_{ad}}^0(u)$ where $u$ is in the
image of $\g_E$ under $H_{r'}@>>>(H_{r'})_{ad}$; this finite group is
well defined up to isomorphism.
If $\ov{Pr}(E)=\ov{\Pr}$ we define
$\ov{Pr}'(E)=\{r'\in\ov{Pr};\ca_{r',E}\cong\ca_{0,E}\}$; this is a
subset of $\ov{Pr}$ with finite complement.

We define a finite collection $c(E)$ of finite groups as follows.

If $\ov{Pr}(E)=\ov{\Pr}=\ov{Pr}'(E)$, then $c(E)$ consists of $\ca_{0,E}$.

If $\ov{Pr}(E)=\ov{\Pr}\ne\ov{Pr}'(E)$, then $c(E)$ consists of 
$\{\ca_{r',E};r'\in\ov{\Pr}-\ov{Pr}'(E)\}$; one can verify that for
$r'\ne r''$ in $\ov{\Pr}-\ov{Pr}'(E)$, we have $\ca_{r',E}\not\cong\ca_{r'',E}$.

If $\ov{Pr}(E)\ne\ov{\Pr}$, then $\ov{Pr}(E)$ consists of a single element
$r'_0\in\ov{Pr}$ (we have necessarily $r'_0\ne0$); then $c(E)$ consists of $\ca_{r'_0,E}$.

If $H$ is a torus, then $c(E)$ consists of $\{1\}$.
If $H_{ad}$ is of type $A,B,C$ or $D$, then $c(E)$ consists of a single
group and this is a product of cyclic groups of order $2$.
If $H_{ad}$ is of exceptional type then $c(E)$ consists of one of the
following groups:

(a) $1,\cc_2,\cc_2\T\cc_2,S_3,\D_8,S_3\T\cc_2,S_5$
\nl
or one of the pair of groups:

(b) $(\cc_2,\cc_3),(\cc_4,\cc_3),(\cc_2\T\cc_2,\cc_2\T\cc_3)$
\nl
or the triple of groups:

(c) $(\cc_4,\cc_3,\cc_5)$.
\nl
(See the tables in \S2.) Here $\cc_m$ denotes a cyclic group of oder $m$,
$S_m$ denotes the symmetric group in $m$ letters, $\D_8$ denotes a
dihedral group of order $8$.

We now define a finite set $c(E)^*$ as follows.
If $c(E)$ consists of a single group $\G$ then $c(E)^*=\hat\G$.
(For a finite group $\G$ we denote by $\hat\G$ the set of isomorphism
classes of irreducible representations of $\G$ over $K_r$.)

If $c(E)$ consists of two groups $\G,\G'$ (see (b)), then $\G''=\ca_{0,E}$ is well
defined and is a quotient of both $\G,\G'$.
(We have $\G''=1,1,S_2$ respectively in the three cases in (b).) Hence we can regard
$\hat\G''$ as a subset of $\hat\G$ and also as a subset of $\hat\G'$. We define
$c(E)^*=(\hat\G-\hat\G'')\sqc(\hat\G'-\hat\G'')\sqc\hat\G''$.

If $c(E)$ consists of three groups (see (c)) we define
$c(E)^*=\sqc_{m\in[1,6]}\hat\cc_m^!$
where
$\hat\cc_m^!$ consists of the faithful irreducible representations of
$\cc_m$.
(This case occurs only when $H_{ad}$ is of type $E_8$ and when $E=1$. The fact that
$\hat\cc_6^!$ enters in the definition should be connected to the fact that $6$
appears as a coefficient of the highest root of $H$.)

The following theorem can be deduced from the definitions using the results in \S2.
\proclaim{Theorem 1.11}Let $r\in\ov{Pr}$. There exists a bijection
$$CS(H_r)@>\si>>\sqc_{E\in\Irr_*(W)}c(E)^*$$
which makes the following diagram commutative:
$$\CD
CS(H_r)@>\si>>\sqc_{E\in\Irr_*(W)}c(E)^*\\
@V\t VV           @VVV             \\
Str(H_r)@>\si>>\Irr_*(W) \endCD$$
(The left vertical map is as in 0.1(a); the right vertical map is the
obvious one; the lower horizontal map is as in 1.6(a).)
\endproclaim

\head 2. Examples\endhead
\subhead 2.1\endsubhead
Assume that $H_{ad}$ is of type $A_{n-1}$, $n\ge2$. We have

$CS'(H_r)=\{(\emp,E',K_r);E'\in\Irr(W)\}$, $\Irr_*(W)=\Irr(W)$. 
\nl
In this case $\un\t_r$ is the bijection  $(\emp,E',K_r)\m E'$.

\subhead 2.2\endsubhead
Assume that $H_{ad}$ is of type $D_n$, $n\ge4$ or $B_n$, $n\ge3$, or
$C_n$, $n\ge2$. If $H_{ad}$ is of type $D_n$, let $CS''(H)$ be the set of pairs $(J,E')$ where 

$J$ is either $\emp$ (so that $N_W(W_J)/W_J=W$)
 or $J$ is such that $W_J$ is of type $D_{4k^2}$ for some 
$k\ge1$ with $4k^2\le n$ (so that $N_W(W_J)/W_J$ is a Weyl group
of type $B_{n-4k^2}$) and $E'\in\Irr(N_W(W_J)/W_J)$.
(We use the convention that a Weyl group of type $B_0$ is $\{1\}$.)

If $H_{ad}$ is of type $B_n$ or $C_n$, let $CS''(H)$ be the set of
pairs $(J,E')$ where 

$J$ is either $\emp$ (so that $N_W(W_J)/W_J=W$)
 or $J$ is such that $W_J$ is of type $B_{k(k+1)}$ for some 
$k\ge1$ with $k(k+1)\le n$ (so that $N_W(W_J)/W_J$ is a Weyl group
of type $B_{n-k(k+1)}$) and $E'\in\Irr(N_W(W_J)/W_J)$.

In any case we have a bijection $CS'(H_r)@>\si>>CS''(H)$ given by
$(J,E',A')\m(J,E')$.
Moreover we have $\Irr_*(W)=\Irr_2(W)$. Hence the map $\un\t_r$ can be viewed 
as a map 
$$CS''(H)@>>>Irr_2(W).\tag a$$
Now $CS''(H)$ can also be viewed as the set of pairs consisting of 
a cuspidal $J$ and a cuspidal local system on a unipotent class in $L_{J,2}$.
The generalized Springer correspondence \cite{L84} attaches to such a pair
a unipotent class in $H_2$ and an irreducible local system on it.

By forgetting this last local system and by identifying $\cu(H_2)$
with $\Irr_2(W)$ 
via $e_2$ (see 1.5), we obtain a map $CS''(H_r)@>>>\Irr_2(W)$ which, on the one hand,
is explicitly computed in \cite{LS85} in terms of certain types of symbols and, 
on the other hand, it coincides with the map (a).

\subhead 2.3\endsubhead
In 2.4-2.8 we describe the map $\un\t_r$ in terms of tables in the case where 
$H_{ad}$ is of type $G_2,F_4,E_6,E_7$ or $E_8$. The tables are computed using
results in \cite{S85} with one indeterminacy in type $E_8$ being removed by \cite{H22}.

In each case the table consists of a sequence of rows. There is one row
for each $E\in\Irr_*(W)$; it is written as $()'.......()''$
where $()'$ represents the fibre of $\un\t_r$ over $E$ and $()''$ is a
sequence of finite groups of which the boxed ones describe $c(E)$.

The elements of $()'$ are written as symbols 
$(J,E',d)_{\sha=n}$. Such a symbol stands for the $n$ triples
$(J,E',A')$ in $CS'(H_r)$ with $J,E'$ fixed
and $A'$ running through the set
$CS^\emp_d(L_{J,r})$ (assumed to have $n\ge1$ elements).
When $n=1$ we omit the subscript $\sha=n$.
We specify $J$ by indicating the type of $W_J$.
(For example, in the table for $E_8$ in 2.8, the row of $8_1$ contains
an item $(E_6,\e_c,0)_{\sha}$ which stands for two objects; in the
triple  $(E_6,\e_c,0)$, $E_6$ represents a subset of type $E_6$ of the
simple reflections, $\e_c$ is a certain representation of a Weyl group of type $G_2$ and $0$ represents the dimension of a certain variety.)
When $J=\emp$ we must have $d=0,n=1$ and we write $E'$ instead of $(J,E',d)$. 
Note that the first entry in $()'$ is $E$ itself.

The groups in $()''$ are as follows.
If $\ov{Pr}(E)=\ov{Pr}=\ov{Pr}'(E)$ then 
$()''$ consists of the single group in $c(E)$ put inside a box.

If $\ov{Pr}(E)=\ov{Pr}\ne\ov{Pr}'(E)$ then 
$()''$ is $\G,\G',(\G'')$ where $\G=\ca_{2,E},\G'=\ca_{3,E},\G''=\ca_{0,E}$;
the boxed entries $\G$ or $\G''$ or both represent the set $c(E)$; an
exception is when $E=1$ in type $E_8$: in this case $()''$ is $\G,\G',\G'',(1)$
where $\G=\ca_{2,E}=\cc_4,\G'=\ca_{3,E}=\cc_3,\G''=\ca_{5,E}=\cc_5,
\ca_{0,E}=1$
and $c(E)$ consists of $\G,\G',\G''$ (all boxed).

If $\ov{Pr}(E)\ne\ov{Pr}$ then $\ov{Pr}-\ov{Pr}(E)=\{r'_0\}$ where
$r'_0\in\{2,3\}$.
If $r'_0=2$ then $()''$ is $\G,-,(-)$ where $\G=\ca_{2,E}$ and $c(E)$
consists of $\G$ (it is boxed);
if $r'_0=3$ then $()''$ is $-,\G',(-)$ where $\G'=\ca_{3,E}$
and $c(E)$ consists of $\G'$ (it is boxed).

If $H_{ad}$ is of type $G_2$, we have $J=\emp$ or $W_J=W$ with
$N_W(W_J)/W_J=\{1\}$.

If $H_{ad}$ is of type $F_4$, we have $J=\emp$ or $W_J$ of type
$B_2$ with $N_W(W_J)/W_J$ of type $B_2$ or $W_J=W$ with
$N_W(W_J)/W_J=\{1\}$.

If $H_{ad}$ is of type $E_6$, we have $J=\emp$ or $W_J$ of type
$D_4$ with $N_W(W_J)/W_J$ of type $A_2$ or $W_J=W$
with $N_W(W_J)/W_J=\{1\}$.

If $H_{ad}$ is of type $E_7$, we have $J=\emp$ or $W_J$ of type
$D_4$ with $N_W(W_J)/W_J$ of type $B_3$ or $W_J$ of type $E_6$
with $N_W(W_J)/W_J$ of type $A_1$ or $W_J=W$
with $N_W(W_J)/W_J=\{1\}$.

If $H_{ad}$ is of type $E_8$, we have $J=\emp$ or $W_J$ of type $D_4$
with $N_W(W_J)/W_J$ of type $F_4$ or $W_J$ of type $E_6$ with
$N_W(W_J)/W_J$ of type $G_2$ or $W_J$ of type $E_7$ with $N_W(W_J)/W_J$
of type $A_1$ or $W_J=W$ with $N_W(W_J)/W_J=\{1\}$.

The notation for the elements of $\Irr(W)$ or $\Irr(N_W(W_J)/W_J)$ is
taken from \cite{S85} with one note of caution. In the case where
$H_{ad}$ is of type $E_8$ and $W_J$ is of type $E_6$, the two
$2$-dimensional irreducible representations $\th',\th''$ of
$N_W(W_J)/W_J$ which appear in the generalized Springer correspondence with $r=3$
are identified in \cite{S85} only up to order. This indeterminacy is
removed in \cite{H22} which shows that $\th'$ is the reflection
representation.

\subhead 2.4. Table for $G_2$\endsubhead

$\e$.......$\bx{1}$

$\e_l$.......$-,\bx{1},(-)$

$\e_c$.......$\bx{1}$

$\th''$.......$\bx{1}$

$\th',(G_2,1,1)$......$S_3,\bx{\cc_2},(S_3)$

$1,(G_2,1,0)_{\sha=3}$......$\bx{\cc_2,\cc_3},(1)$

\subhead 2.5. Table for $F_4$\endsubhead

$\c_{1,4}$.......$\bx{1}$

$\c_{2,4}$.......$\bx{1}$

$\c_{2,2}$.......$\bx{1},-,(-)$

$\c_{4,4}$......$\bx{1},\cc_2,(\cc_2)$

$\c_{9,4}$.......$\bx{1}$

$\c_{8,4},\c_{1,2}$.......$\bx{\cc_2}$

$\c_{8,2},\c_{1,3}$...... $\bx{\cc_2},1,(1)$

$\c_4,(B_2,\e,0)$........$\bx{\cc_2},-,(-)$

$\c_{4,3}$.......$\bx{1},-,(-)$

$\c_{4,2}$.......$\bx{1}$

$\c_{9,3}$.......$\bx{1},-,(-)$

$\c_{9,2}$.......$\bx{1},\cc_2,(\cc_2)$

$\c_{6,1}$.......$\bx{1}$

$\c_{16}$.......$\bx{1},\cc_2,(\cc_2)$

$\c_{12},\c_{6,2}, (F_4,1,4)$.......$\bx{S_3},S_4,(S_4)$

$\c_{8,3},(B_2,\e_l,0)$.......$\bx{\cc_2},1,(1)$

$\c_{8,1},(B_2,\e_c,0)$.......$\bx{\cc_2},1,(1)$

$\c_{9,1},\c_{2,1},\c_{2,3},(B_2,\th,0),(F_4,1,2)$.......
$\bx{\D_8},\cc_2,(\cc_2)$

$\c_{4,1},(F_4,1,1)$......$\bx{\cc_2}$

$\c_{1,1},(B_2,1,0),(F_4,1,0)_{\sha=4}$.......$\bx{\cc_4,\cc_3},(1)$

\subhead 2.6. Table for $E_6$\endsubhead

$1_{36}$ .......$\bx{1}$

$6_{25}$.......$\bx{1}$

$20_{20}$.......$\bx{1}$

$15_{16}$.......$\bx{1}$

$30_{15},15_{17}$.......$\bx{\cc_2}$

$64_{13}$.......$\bx{1}$

$24_{12}$.......$\bx{1}$

$60_{11}$.......$\bx{1}$

$81_{10}$.......$\bx{1}$

$10_9$.......$\bx{1}$

$60_8$.......$\bx{1}$

$80_7,90_8,20_{10}$ .......$\bx{S_3}$

$81_6$.......$\bx{1}$

$24_6,(D_4,\e,0)$.......$\bx{S_2},1,(1)$

$60_5$.......$\bx{1}$

$64_4$.......$\bx{1}$

$15_4$.......$\bx{1}$

$30_3,15_5$.......$\bx{\cc_2}$ 

$20_2,(D_4,\ph,0)$.......$\bx{\cc_2},1,(1)$

$6_1$.......$\bx{1}$

$1_0,(D_4,1,0),(E_6,1,0)_{\sha=2}$.......$\bx{\cc_2,\cc_3},(1)$

\subhead 2.7. Table for $E_7$\endsubhead

$1_{63}$.......$\bx{1}$

$7_{46}$.......$\bx{1}$

$27_{37}$.......$\bx{1}$

$21_{36}$.......$\bx{1}$

$35_{31}$.......$\bx{1}$

$56_{30},21_{33}$.......$\bx{\cc_2}$

$15_{28}$.......$\bx{1}$

$120_{25},105_{28}$.......$\bx{\cc_2}$

$189_{22}$.......$\bx{1}$

$105_{21}$.......$\bx{1}$

$168_{21}$.......$\bx{1}$

$210_{21}$.......$\bx{1}$

$189_{20}$.......$\bx{1}$

$70_{18}$.......$\bx{1}$

$280_{17}$.......$\bx{1}$

$315_{16},280_{18},35_{22}$ .......$\bx{S_3}$

$216_{16}$.......$\bx{1}$

$405_{15},189_{17}$.......$\bx{\cc_2}$

$105_{15},(D_4,(0,1^3),0)$.......$\bx{\cc_2},1,(1)$

$84_{15}$.......$\bx{1},-,(-)$

$378_{14}$.......$\bx{1},\cc_2,(\cc_2)$ 

$210_{13}$.......$\bx{1}$

$420_{13},336_{14}$.......$\bx{\cc_2}$ 

$84_{12},(D_4,(1^3,0),0)$.......$\bx{\cc_2}$

$105_{12}$.......$\bx{1}$

$512_{11},512_{12}$.......$\bx{\cc_2}$

$210_{10}$.......$\bx{1}$

$420_{10},336_{11}$.......$\bx{\cc_2}$

$378_9$.......$\bx{1}$

$216_9$.......$\bx{1}$

$70_9$.......$\bx{1}$

$280_8$.......$\bx{1}$

$405_8,189_{10}$.......$\bx{\cc_2}$ 

$189_7,(D_4,(1,1^2),0)$.......$\bx{\cc_2},1,(1)$

$315_7,280_9,35_{13}$.......$\bx{S_3}$

$168_6,(D_4,(1^2,1),0)$.......$\bx{\cc_2},1,(1)$

$210_6,(D_4,(0,21),0)$.......$\bx{\cc_2},1,(1)$

$105_6,15_7$.......$\bx{\cc_2},1,(1)$

$189_5$.......$\bx{1},\cc_2,(\cc_2)$

$35_4,(D_4,(21,0),0)$.......$\bx{\cc_2},1,(1)$

$120_4,105_5$.......$\bx{\cc_2}$ 

$21_3,(D_4,(1,2),0),(E_6,1,0)_{\sha=2}$.......$\bx{\cc_2,\cc_3},(1)$

$56_3,21_6$.......$\bx{\cc_2}$ 

$27_2,(D_4,(2,1),0)$.......$\bx{\cc_2},1,(1)$ 

$7_1,(D_4,(0,3),0)$.......$\bx{\cc_2},1,(1)$  

$1_0,(D_4,(3,0),0),(E_6,1,0)_{\sha=2},(E_7,1,0)_{\sha=2}$
.......$\bx{\cc_4,\cc_3},(1)$

\subhead 2.8. Table for $E_8$\endsubhead

$1_{120}$.......$\bx{1}$

$8_{91}$.......$\bx{1}$

$35_{74}$.......$\bx{1}$

$84_{64}$.......$\bx{1}$

$112_{63},28_{68}$.......$\bx{\cc_2}$
 
$50_{56}$.......$\bx{1}$

$210_{52},160_{55}$.......$\bx{\cc_2}$
 
$560_{47}$.......$\bx{1}$

$567_{46}$.......$\bx{1}$

$400_{43}$.......$\bx{1}$

$700_{42},300_{44}$.......$\bx{\cc_2}$  
 
$448_{39}$.......$\bx{1}$ 

$1344_{38}$.......$\bx{1}$

$1400_{37},1008_{39},56_{49}$.......$\bx{S_3}$
 
$175_{36}$.......$\bx{1}$

$525_{36},(D_4,\c_{1,4},0)$.......$\bx{\cc_2},1,(1)$
 
$1050_{34}$.......$\bx{1}$

$1400_{32},1575_{34},350_{38}$.......$\bx{S_3}$

$972_{32}$.......$\bx{1},-,(-)$

$3240_{31}$.......$\bx{1},\cc_2,(\cc_2)$    
 
$2268_{30},1296_{33}$.......$\bx{\cc_2}$
 
$1400_{29}$.......$\bx{1}$

$2240_{28},840_{31}$.......$\bx{\cc_2}$
 
$700_{28}, (D_4,\c_{2,2},0)$.......$\bx{\cc_2},1,(1)$

$840_{26}$.......$\bx{1}$

$4096_{26},4096_{27}$.......$\bx{\cc_2}$

$2800_{25},2100_{28}$.......$\bx{\cc_2}$
 
$4200_{24},3360_{25}$.......$\bx{\cc_2}$
 
$168_{24},(D_4,\c_{1,3},0)$.......$\bx{\cc_2},-,(-)$

$4536_{23}$.......$\bx{1}$

$2835_{22}$.......$\bx{1}$

$6075_{22}$.......$\bx{1}$

$3200_{32}$.......$\bx{1}$

$4200_{21}$.......$\bx{1},\cc_2,(\cc_2)$
 
$5600_{21},2400_{23}$.......$\bx{\cc_2}$

$420_{20}$.......$\bx{1}$

$2100_{20},(D_4,\c_{4,4},0)$.......$\bx{\cc_2},1,(1)$ 
 
$1344_{19}$.......$\bx{1}$

$2016_{19}$.......$\bx{1}$

$3150_{18},1134_{20}$.......$\bx{\cc_2}$
 
$4200_{18},2688_{20}$.......$\bx{\cc_2}$
 
$7168_{17},5600_{19},448_{25}$.......$\bx{S_3}$
 
$3200_{16},(D_4,\c_{8,2},0)$.......$\bx{\cc_2},1,(1)$

$4480_{16},5670_{18},4536_{18},1400_{20},1680_{22},70_{32},(E_8,1,16)$.......$\bx{S_5}$ 

$5600_{15},2400_{17},(D_4,\c_{9,4},0),(D_4,\c_{2,4},0)$.......
$\bx{\cc_2\T\cc_2},\cc_2,(\cc_2)$

$4200_{15},700_{16}$.......$\bx{\cc_2},1,(1)$

$2835_{14}$.......$\bx{1}$

$6075_{14}$.......$\bx{1},\cc_2,(\cc_2)$

$840_{14},(D_4,\c_{4,3},0)$ .......$\bx{\cc_2},-,(-)$

$4536_{13}$.......$\bx{1},\cc_2,(\cc_2)$

$2800_{13},2100_{16}$ .......$\bx{\cc_2}$ 

$972_{12},(D_4,\c_{9,3},0)$.......$\bx{\cc_2},1,(1)$ 

$4200_{12},3360_{13}$.......$\bx{\cc_2}$

$525_{12},(D_4,\c_{8,4},0),(E_6,\e,0)_{\sha=2}$.......$\bx{\cc_2,\cc_3},(1)$

$175_{12}$.......$-,\bx{1},(-)$

$1400_{11}$.......$\bx{1}$

$4096_{11},4096_{12}$.......$\bx{\cc_2}$ 

$2268_{10},1296_{13}$.......$\bx{\cc_2}$

$2240_{10},840_{13}$.......$S_3,\bx{\cc_2},(S_3)$ 

$1050_{10},(D_4,\c_4,0)$.......$\bx{\cc_2},-,(-)$ 

$3240_9$.......$\bx{1},\cc_2,(\cc_2)$ 

$448_9,(D_4,\c_{6,1},0),(E_6,\e_l,0)_{\sha=2}$.......$\bx{\cc_2,\cc_3},(1)$

$1344_8,(D_4,\c_{16},0)$.......$\bx{\cc_2},1,(1)$

$1400_8,1575_{10},350_{14}$.......$\bx{S_3}$ 

$1400_7,1008_9,56_{19},(D_4,\c_{12},0),(D_4,\c_{6,2},0),(E_8,1,7)$
.......$\bx{S_3\T\cc_2},S_3,(S_3)$ 

$400_7,(D_4,\c_{2,3},0)$.......$\bx{\cc_2},1,(1)$

$700_6,300_8,50_9,(D_4,\c_{8,3},0),(E_8,1,6)$.......$\bx{\D_8},\cc_2,
(\cc_2)$
 
$567_6,(D_4,\c_{9,2},0)$.......$\bx{\cc_2},1,(1)$ 

$560_5,(D_4,\c_{4,2},0)$.......$\bx{\cc_2}$ 

$210_4,160_7$.......$\bx{\cc_2}$ 

$84_4,(D_4,\c_{9,1},0),(E_6,\th'',0)_{\sha=2},(E_7,1,0)_{\sha=2}$
.......$\bx{\cc_4,\cc_3},(1)$

$112_3,28_8,(D_4,\c_{8,1},0),(D_4,\c_{1,2},0),(E_6,\th',0)_{\sha=2},
(E_8,1,3)_{\sha=2}$

.......$\bx{\cc_2\T\cc_2,\cc_2\T\cc_3},(\cc_2)$

$35_2,(D_4,\c_{4,1},0)$.......$\bx{\cc_2},1,(1)$

$8_1,(D_4,\c_{2,1},0),(E_6,\e_c,0)_{\sha=2},(E_8,1,1)_{\sha=2}$
.......$\bx{\cc_4,\cc_3},(1)$

$1_0,(D_4,\c_{1,1},0), (E_6,1,0)_{\sha=2},
(E_7,1,0)_{\sha=2},(E_8,1,0)_{\sha=6}$.......$\bx{\cc_4,\cc_3,\cc_5},(1)$ 
\head 3. Complements\endhead
\subhead 3.1\endsubhead
The restriction of the map $\t$ in 0.1(a) to $CS^\emp(H_r)$ has an
alternative definition. Namely, for $A\in CS^\emp(H_r)$, there is a
unique stratum $X\in Str(H_r)$ such that $\s_A\sub X$ (notation of 1.1);
we have $\t(A)=X$.

\subhead 3.2\endsubhead
The results in this subsection can be used to verify 1.7.
In the examples below we assume that $H$ is semisimple, $H\ne\{1\}$,
and for $A\in CS^\emp(H_r)$ we denote by $s$ the semisimple part of an 
element of $\s_A$. We describe the structure
of $Z_{H_r}^0(s)$ in various cases. We also specify the value of $X$
in 1.7.

If $H$ is of type $C_n$ with $n=k(k+1),k\ge1$ then:

if $r\ne 2$ then $Z_{H_r}^0(s)$ is of type $C_{n/2}\T C_{n/2}$;
if $r=2$ then $Z_{H_r}^0(s)=H_r$. Thus $X$ is as in 1.7(i).

If $H$ is of type $B_n$ with $n=k(k+1),k\ge1$ then:

if $r\ne2$ then $Z_{H_r}^0(s)$ is of type $B_a\T D_b$
where $(2a+1,2b)=((k+1)^2,k^2)$ if $k$ is even
and $(2a+1,2b)=(k^2,(k+1)^2)$ if $k$ is odd; if $r=2$ then
$Z_{H_r}^0(s)=H_r$. Thus $X$ is as in 1.7(i).

If $H$ is of type $D_n$ with $n=4k^2,k\ge1$ then:

if $r\ne2$ then $Z_{H_r}^0(s)$ is of type $D_{2k^2}\T D_{2k^2}$; if
$r=2$ then $Z_{H_r}^0(s)=H_r$. Thus $X$ is as in 1.7(i).

If $H$ is of type $G_2$ and $A\in CS^\emp_d(H_r)$ then:

if $d=1$ then $Z_{H_r}^0(s)=H_r$ (thus $X$ is as in 1.7(ii));
if $d=0$ and $r\n\{2,3\}$ then $Z_{H_r}^0(s)$ is of type $A_2$ for two
values of $A$ and of type $A_1\T A_1$ for the third value of $A$;
if $d=0$ and $r=2$ then $Z_{H_r}^0(s)$ is of type $A_2$ for two
values of $A$ and is $H_r$ for the third value of $A$;
if $d=0$ and $r=3$ then $Z_{H_r}^0(s)$ is $H_r$ for two
values of $A$ and is of type $A_1\T A_1$ for the third value of $A$.
Thus $X$ is as in 1.7(iv).

If $H$ is of type $F_4$ and $A\in CS^\emp_d(H_r)$ then:

if $d=4$ then $Z_{H_r}^0(s)=H_r$ (thus $X$ is as in 1.7(ii));

if $d=2$, $r\ne2$,  then $Z_{H_r}^0(s)$ is of type $B_4$;
if $d=2$, $r=2$,  then $Z_{H_r}^0(s)=H_r$ (thus $X$ is as in 1.7(i));

if $d=1$, $r\ne2$, then $Z_{H_r}^0(s)$ is of type $C_3\T A_1$;
if $d=1$, $r=2$,  then $Z_{H_r}^0(s)=H_r$ (thus $X$ is as in 1.7(i)); 

if $d=0$, $r\n\{2,3\}$, then $Z_{H_r}^0(s)$ is of type $A_2\T A_2$
for two values of $A$ and of type $A_3\T A_1$ for the other two values
of $A$; if $d=0$, $r=2$, then $Z_{H_r}^0(s)$ is of type $A_2\T A_2$
for two values of $A$ and is $H_r$ for the other two values of $A$;
if $d=0$, $r=3$, then $Z_{H_r}^0(s)$ is of type $A_3\T A_1$
for two values of $A$ and is $H_r$ for the other two values of $A$.
Thus $X$ is as in 1.7(iv).

If $H$ is of type $E_6$ then:

if $r\ne3$, then $Z_{H_r}^0(s)$ is of type $A_2\T A_2\T A_2$;
if $r=3$,  then $Z_{H_r}^0(s)=H_r$. Thus $X$ is as in 1.7(i).

If $H$ is of type $E_7$ then $d=0$ and:

if $r\ne2$, then $Z_{H_r}^0(s)$ is of type $A_3\T A_3\T A_1$;
if $r=2$, then $Z_{H_r}^0(s)=H_r$. Thus $X$ is as in 1.7(i).

If $H$ is of type $E_8$ and $A\in CS^\emp_d(H_r)$ then:

if $d=16$ then $Z_{H_r}^0(s)=H_r$ (thus $X$ is as in 1.7(ii)).

if $d=7$ and $r\ne2$ then $Z_{H_r}^0(s)$ is of type $E_7\T A_1$;
if $d=7$ and $r=2$ then $Z_{H_r}^0(s)=H_r$ (thus $X$ is as in 1.7(i));

if $d=6$ and $r\ne2$ then $Z_{H_r}^0(s)$ is of type $D_8$;
if $d=6$ and $r=2$ then $Z_{H_r}^0(s)=H_r$ (thus $X$ is as in 1.7(i));

if $d=3$ and $r\ne3$ then $Z_{H_r}^0(s)$ is of type $E_6\T A_2$;
if $d=3$ and $r=3$ then $Z_{H_r}^0(s)=H_r$ (thus $X$ is as in 1.7(i));

if $d=1$ and $r\ne2$ then $Z_{H_r}^0(s)$ is of type $D_5\T A_3$;
if $d=1$ and $r=2$ then $Z_{H_r}^0(s)=H_r$ (thus $X$ is as in 1.7(i));

if $d=0$ and $r\n\{2,3,5\}$ then
$Z_{H_r}^0(s)$ is of type $A_4\T A_4$ for four values of $A$ and of type
$A_5\T A_2\T A_1$ for two values of $A$; if $d=0$ and $r=5$ then 
$Z_{H_r}^0(s)$ is $H_r$ for four values of $A$ and of type
$A_5\T A_2\T A_1$ for two values of $A$;
if $d=0$ and $r=3$ then $Z_{H_r}^0(s)$ is of type $A_4\T A_4$ for 
four values of $A$ and of type $E_7\T A_1$ for two values of $A$;
if $d=0$ and $r=2$ then $Z_{H_r}^0(s)$ is
of type $A_4\T A_4$ for four values of $A$ and of type
$E_6\T A_2$ for two values of $A$. Thus $X$ is as in 1.7(iv).

The results in this section (for $H$ of type $E_8$ and $d=0$)
contradict the statement (f) on p.351 in \cite{Sh95}. (Indeed, if $r=2$
there is no semisimple $s\in H_2$ with $Z_{H_2}(s)$ of type
$A_5\T A_2\T A_1$.)

\subhead 3.3\endsubhead
Let $H^*$ be a connected reductive group over $\CC$ of type dual to that
of $H$. 
Let $\CC\SS(H)$ be the set of isomorphism classes of (not
necessarily unipotent) character sheaves on $H$.
We now assume that $Z_H=Z_H^0$. It is known that we can identify
$\CC\SS(H)$ with $\sqc_s CS(Z_{H^*}(s)^*)$ where $s$ runs over
the semisimple elements of finite order of $H^*$ up to conjugacy.
Using 0.1(a), we obtain a surjective map
$\CC\SS(H)@>>>\sqc_sStr(Z_{H^*}(s)^*)$.
From \cite{L15} we can identify $Str(Z_{H^*}(s)^*)=Str(Z_{H^*}(s))$.
Hence we obtain a surjective map
$$\CC\SS(H)@>>>\sqc_sStr(Z_{H^*}(s)).$$

\subhead 3.4\endsubhead
Let $\cx'_H$ be the set of numbers which appear as coefficients of the
highest root of $H$; let $\cx_H=\cx'_H\cup\{1\}$. Note that $\cx_H$
consists
of the numbers $1,2,\do,z_H$ where $z_H=1$ for $H$ of type $A$, $z_H=2$
for $H$ of type $B,C,D$, $z_H=3$ for $H$ of type $G_2$ or $E_6$, $z_H=4$
for $H$ of type $F_4$ or $E_7$, $z_H=6$ for $H$ of type $E_8$. 

We note the following property (which can be verified from the results
in \S2).

(a) {\it The fibre of $\t$ in 0.1(a) at the stratum
consisting of regular
elements in $H_r$ (or equivalently at $E=1\in\Irr_*(W)$) is in bijection
with the set
$$\sqc_{m\in\NN;1\le m\le z_H}\r_m\tag b$$
where $\r_m$ is the set of primitive $m$-th roots of $1$ in $K_r$.}
\nl
It is remarkable that the set (b) appears also in a quite different
situation. Let $r,q,H_r(F_q)$ be as in 0.1(b). We can view $H_r(F_q)$
as a fixed point set of a Frobenius map $F:H_r@>>>H_r$. For any
$w\in W$ let $X_w$ be the variety attached to $H_r,F,w$ in \cite{DL76}.
Now $F$ acts on the cohomology with compact support $H^i_c(X_w)$ of
$X_w$ and in particular on $H^{|w|}_c(X_w)$. (We denote by $|w|$ the
length of $w$.) Let $w$ be a Coxeter
element of minimal length in $W$. From \cite{L76} it is known that the
$F$-action on $H^{|w|}_c(X_w)$ is semisimple and that the eigenspaces
are irreducible (unipotent) representations of $H_r(F_q)$. These
unipotent representations are in bijection with the character
sheaves in the fibre of $\t$ at $E=1$.
(We use the usual bijection $Un(H_r(F_q))\lra CS(H_r)$
composed with the involution of $Un(H_r(F_q))$ which interchanges
``small'' representations with ``big'' representations.)
The eigenvalues of the
$F$-action are listed in \cite{L76, p.146,147}. It turns out that

(c) these eigenvalues are exactly the roots of $1$ in (b) times
integral powers of $q^{1/2}$.

\subhead 3.5\endsubhead
Let $cl(W)$ be the set of conjugacy classes in $W$. In \cite{L15,\S4},
a surjective map $\Ph:cl(W)@>>>\Irr_*(W)$ is defined. In
\cite{L15, 4.10} a map $\Irr_*(W)@>>>cl(W)$, $E@>>>C_E$, is described;
it is such that $\Ph(C_E)=E$ for all $E\in\Irr_*(W)$, hence its image
$cl_*(W)\sub cl(W)$ is such that $\Ph$ restricts to a bijection
$cl_*(W)@>\si>>\Irr_*(W)$.
This allows us to identify the sets $cl_*(W),\Irr_*(W)=Str(H_r)$,
so that $\t:CS(H_r)@>>>Str(H_r)$ becomes a surjective map

(a) $\t':Un(H_r(F_q))@>>>cl_*(W)$
\nl
(with $r,q,H_r(F_q)$ as in 0.1(b), see 0.1(c)).
Let $Un^\emp(H_r(F_q))$ be the subset of $Un(H_r(F_q))$ consisting of
unipotent cuspidal representations.

One can verify that the restriction of $\t'$ to $Un^\emp(H_r(F_q))$
coincides with the map $\r\m C_\r$ in \cite{L02, 2.17}. From
\cite{L02} we see that

(b) for any $\r\in Un^\emp(H_r(F_q))$, $\r$ appears with multiplicity
$1$ in $H^{|w|}_c(X_w)$ (notation of 3.4) where $w$ is an element of
minimal length in $\t'(\r)$.

\subhead 3.6\endsubhead
In this subsection we assume that $H_r$ is the symplectic group
with $W$ of type $B_2$. The simple reflections $s_1,s_2$ satisfy
$s_1s_2s_1s_2=s_2s_1s_2s_1$. We have $\Irr_*(W)=\Irr(W)=\Irr_2(W)$;
this set consists of $1,\r,\e_1,\e_2,\e$ where $\r$ is the reflection
representation, $\e$ is the sign representation and $\e_1,\e_2$ are the
one dimensional representations other than $1,\e$. We have
$cl(W)=cl_*(W)$. (The numbering of $s_1,s_2$ and of $\e_1,\e_2$
is chosen so that (a),(b) below hold.)

The conjugacy classes in $W$ are 
$(1),(s_1),(s_2),(s_1s_2),(s_1s_2s_1s_2)$ where $(w)$ is the conjugacy
class of $w\in W$. The bijection $\Irr_*(W)@>>>cl_*(W)$ is given by

(a) $1\m(s_1s_2)$, $\r\m(s_1s_2s_1s_2)$, $\e_1\m(s_1)$, $\e_2\m(s_2)$,
$\e\m(1)$.
\nl
There are five strata; they are indexed by the elements of $\Irr_*(W)$;
we denote them by $\s(1),\s(\r),\s(\e_1),\s(\e_2),\s(\e)$. Here

(b) $\s(1)$ is the union of all conjugacy classes of dimension $8$;
$\s(\r)$ is the union of all conjugacy classes of dimension $6$;
$\s(\e_1)$ is a conjugacy class of dimension $4$ (a semisimple one if
$r\ne2$ and a unipotent one if $r=2$; $\s(\e_2)$ is a union of one (if
$r=2$) or two (if $r\ne2$) conjugacy classes of dimension $4$; $\s(\e)$
is the centre of $H_r$.

$CS(H_r)$ consists of six objects: $\bbq[1],\bbq[\r],\bbq[\e_1],
\bbq[\e_2],\bbq[\e]$ and $A$ (an object of $CS^\emp(H_r)$).
The map $\t$ is

$\bbq[1]\m1,\bbq[\r]\m\r,\bbq[\e_1]\m\e_1,
\bbq[\e_2]\m\e_2,\bbq[\e]\m\e$, $A\m1$.

\enddocument